\pdfoutput=1
\RequirePackage{ifpdf}
\ifpdf 
\documentclass[pdftex]{sigma}
\else
\documentclass{sigma}
\fi

\usepackage[all]{xy}

\newcommand{\M}{\ensuremath{\mathcal M}}

\newcommand{\CO}{\ensuremath{\mathcal O}}

\newcommand{\cA}{\ensuremath{\mathcal A}}

\newcommand{\GC}{\ensuremath{G_{\mathbb{C}}}}
\newcommand{\MGC}{\ensuremath{\mathcal{M}_{G_{\mathbb{C}}}}}
\newcommand{\MG}{\ensuremath{\mathcal{M}_{G}}}

\newcommand{\C}{\ensuremath{\mathbb C}}

\newcommand{\Z}{\ensuremath{\mathbb Z}}

\newcommand{\e}{{\mathrm e}}
\newcommand{\de}{\mathrm{d}}

\newcommand{\eps}{\epsilon}

\newcommand{\delbar}{\overline{\partial}}

\numberwithin{equation}{section}

\newtheorem{Theorem}{Theorem}[section]

\newtheorem{open}[Theorem]{Open Question}
\newtheorem{Conjecture}[Theorem]{Conjecture}
 { \theoremstyle{definition}
\newtheorem{Definition}[Theorem]{Definition}

 }

\begin{document}

\newcommand{\arXivNumber}{1710.08453}

\renewcommand{\PaperNumber}{037}

\FirstPageHeading

\ShortArticleName{Singular Geometry and Higgs Bundles in String Theory}

\ArticleName{Singular Geometry and Higgs Bundles\\ in String Theory}

\Author{Lara B.~ANDERSON~$^{\dag^1}$, Mboyo ESOLE~$^{\dag^2}$, Laura FREDRICKSON~$^{\dag^3}$\newline and Laura P.~SCHAPOSNIK~$^{\dag^4\dag^5}$}

\AuthorNameForHeading{L.B.~Anderson, M.~Esole, L.~Fredrickson and L.P.~Schaposnik}

\Address{$^{\dag^1}$~Department of Physics and Department of Mathematics, Virginia Tech,\\
\hphantom{$^{\dag^1}$}~Blacksburg, VA 24061, USA} 
\EmailDD{\href{mailto:lara.anderson@vt.edu}{lara.anderson@vt.edu}}

\Address{$^{\dag^2}$~Department of Mathematics, Northeastern University, Boston, MA 02115, USA}
\EmailDD{\href{mailto:j.esole@northeastern.edu}{j.esole@northeastern.edu}}

\Address{$^{\dag^3}$~Department of Mathematics, Stanford University, Stanford, CA 94305, USA}
\EmailDD{\href{mailto:lfredrickson@stanford.edu}{lfredrickson@stanford.edu}}

\Address{$^{\dag^4}$~Department of Mathematics, University of Illinois at Chicago, 60607 Chicago, USA}
\EmailDD{\href{mailto:schapos@uic.edu}{schapos@uic.edu}}
\Address{$^{\dag^5}$~Department of Mathematics, FU Berlin, 14195 Berlin, Germany}

\ArticleDates{Received November 22, 2017, in final form April 13, 2018; Published online April 18, 2018}

\Abstract{This brief survey aims to set the stage and summarize some of the ideas under discussion at the {\it Workshop on Singular Geometry and Higgs Bundles in String Theory}, to be held at the American Institute of Mathematics from October 30th to November 3rd, 2017. One of the most interesting aspects of the duality revolution in string theory is the understanding that gauge fields and matter representations can be described by intersection of branes. Since gauge theory is at the heart of our description of physical interactions, it has opened the door to the geometric engineering of many physical systems, and in particular those involving Higgs bundles. This note presents a curated overview of some current advances and open problems in the area, with no intention of being a complete review of the whole subject.}


\Classification{14D20; 14D21; 53C07; 14H70; 14P25}

\section{Introduction}

One of the most interesting aspects of the duality revolution in string theory is the understanding that gauge fields and matter representations can be described by the intersection of branes. Since gauge theory is at the heart of our description of physical interactions, it has opened the door to the geometric engineering of many physical systems, and in particular those arising from Higgs bundles, whose moduli spaces have become a source of many interesting branes.

In an effort to consolidate and disseminate the variety of different techniques, heuristics, and approaches that have been applied to the study of Higgs bundles and spectral data in recent years by the mathematics and physics communities, we present here a short survey on these subjects, as well as a collection of open problems and ideas revolving around them. This note focuses on two interrelated themes concerning Higgs bundles and the Hitchin fibration, and their interactions with mathematical physics:

\begin{itemize}\itemsep=0pt
\item[(A)] {\it Higgs bundles and algebraic geometry}: invariants of singular spaces, and in particular singular fibers of the Hitchin fibration (Section~\ref{Higgs}); the effect of these fibers on the geometry of the moduli spaces of Higgs bundles, including limits within the Hitchin fibration (Section~\ref{sec-limmits}), and the appearance of Higgs bundles on singular curves (Section~\ref{sec-curves}).

\item[(B)] {\it Hitchin systems and T-branes}: the study of the moduli space of Higgs bundles and its branes through the Hitchin fibration (Section~\ref{sec-branes}), and their appearance within the broader setting of string/F-theory (Section~\ref{sec-string}), and Calabi--Yau elliptic fibrations (Section~\ref{Mboyo_sec}).
\end{itemize}

Although these two themes are closely related, correspondences between them are just in their infancy. In particular, obtaining a global understanding of Higgs bundles over singular curves, and of Higgs bundles which have singular spectral data, would be most beneficial from the perspective of F-theory and superconformal theories in diverse dimensions. We hope these notes will help to further clarify the role that spectral curves and spectral data play in string theory, both for those studying Higgs bundles on Riemann surfaces and for those studying Higgs bundles on higher dimensional spaces.

This short survey is not intended to be a complete overview of the research done in the area, but rather a concise description of certain particular paths of research that are currently receiving much attention, and that present open problems that could be tackled by researchers in different areas of mathematics and physics.

\section{Higgs bundles and the Hitchin fibration}\label{Higgs}

Throughout the paper we shall consider a compact Riemann surface $\Sigma$ of genus $g\geq 2$ with canonical bundle $K=T^*\Sigma$. In what follows, we recall some of the main properties of complex and real Higgs bundles, as well as the associated Hitchin fibration whose structure groups are real or complex subgroups of ${\rm GL}(n,\mathbb{C})$.

 \subsection{Higgs bundles}
We begin by briefly reviewing the notions of Higgs bundles for real and complex groups which are relevant to this paper. Further details can be found in standard references such as Hitchin~\cite{hitchin87,hitchin87b} and Simpson \cite{simpsonhiggs,simpson94,simpson94+}.
Recall that ${\rm GL}(n,\C)$-Higgs bundles of degree 0 on $\Sigma$ are pairs $(E,\Phi)$ where
\begin{itemize}\itemsep=0pt
 \item $E \rightarrow \Sigma$ is a holomorphic vector bundle of rank $n$ and degree $0$,
 \item the Higgs field $\Phi\colon E\rightarrow E \otimes K$, is a holomorphic $K$-valued endomorphism.
 \end{itemize}
 By the work of Hitchin and Simpson, given a polystable Higgs bundle, there is a unique hermitian metric $h$ on $E$,
known as the harmonic metric, solving the so-called {\it Hitchin equations}:
\begin{gather*}
F_{D(\delbar_E, h)}+ [\Phi,\Phi^{*_h}]=0, \qquad \overline{\partial}_{E}\Phi=0,
\end{gather*}
 where $D(\delbar_E{,} h)$ is the Chern connection, i.e., the unique $h$-unitary connection such that \smash{$D^{0,1}{=}\delbar_E$}, the curvature of the Chern connection is denoted by $F_{D(\delbar_E, h)}$, and $\Phi^{*_h}$ represents the hermitian adjoint of~$\Phi$ with respect to the hermitian metric~$h$. The correspondence between pairs $(E, \Phi)$ and triples $(E, \Phi, h)$ is known as the nonabelian Hodge correspondence.
More generally, for a~complex reductive Lie group~$\GC$, one has the following~\cite{hitchin87b}:

\begin{Definition} \label{def1} A $\GC$-Higgs bundle is a pair $(P,\Phi)$, where $P$ is a holomorphic principal $\GC$-bundle, and $\Phi$ is a holomorphic section of ${\rm ad}(P)\otimes K$, where ${\rm ad}(P)$ is the adjoint bundle of~$P$.
\end{Definition}

In this setting, there is a similar nonabelian Hodge correspondence, where the notion of a~hermitian metric is replaced by a reduction of structure of $P$ to the maximal compact subgroup of $\GC$. By considering appropriate stability conditions, one may define the Hitchin moduli space~$\M_{\GC}$ of isomorphism classes of polystable $\GC$-Higgs bundles, which was introduced by Hitchin in~\cite{hitchin87}. Up to gauge equivalence, the points of the moduli space $\M_{\GC}$ represent polystable $\GC$-Higgs bundles on $\Sigma$. Moreover, the through nonabelian Hodge correspondence, points of the moduli space represent solutions of the $\GC$-Hitchin's equations.

Given a real form $G$ of the complex reductive lie group $\GC$, we may define $G$-Higgs bundles as follows. Let $H$ be the maximal compact subgroup of $G$ and consider the Cartan decomposition $\mathfrak{g} = \mathfrak{h}\oplus \mathfrak{m}$ of $\mathfrak{g}$, where $\mathfrak{h}$ is the Lie algebra of $H$, and $\mathfrak{m}$ its orthogonal complement. This induces a~decomposition of the Lie algebra $\mathfrak{g}_\C=\mathfrak{h}^{\mathbb{C}}\oplus \mathfrak{m}^{\mathbb{C}}$ of $G_\C$. Note that the Lie algebras satisfy $ [\mathfrak{h}, \mathfrak{h}]\subset\mathfrak{h}$, $[\mathfrak{h,\mathfrak{m}}]\subset\mathfrak{m}$, $[\mathfrak{m},\mathfrak{m}]\subset \mathfrak{h},$ and there is an induced isotropy representation ${\rm Ad}|_{H^{\mathbb{C}}}\colon H^{\mathbb{C}}\rightarrow {\rm GL}(\mathfrak{m}^{\mathbb{C}})$.

\begin{Definition}\label{def2} A principal $G$-Higgs bundle is a pair $(P,\Phi)$ where
\begin{itemize}\itemsep=0pt
 \item $P$ is a holomorphic principal $H^{\mathbb{C}}$-bundle on $\Sigma$,
 \item $\Phi$ is a holomorphic section of $P\times_{\mathrm{Ad}}\mathfrak{m}^{\mathbb{C}}\otimes K$.
\end{itemize}
\end{Definition}

Similarly to the case of Higgs bundles for complex groups, there are notions of stability, semistability and polystability for $G$-Higgs bundles. One can see that the polystability of a~$G$-Higgs bundle for a group $G\subset {\rm GL}(n,\mathbb{C})$ is equivalent to the polystability of the corresponding ${\rm GL}(n,\mathbb{C})$-Higgs bundle. However, it should be noted that a $G$-Higgs bundle can be stable as a~$G$-Higgs bundle but not as a~${\rm GL}(n,\mathbb{C})$-Higgs bundle. The moduli space of polystable $G$-Higgs bundles on the compact Riemann surface $\Sigma$ shall be denoted by~$\M_G$.

\subsection{The Hitchin fibration} The moduli space $\mathcal{M}_{\GC}$ of $ \GC$-Higgs bundles admits a natural complete hyperk\"ahler metric over its smooth points, and a way of studying it is through the {\it Hitchin fibration}~\cite{hitchin87b}. This fibration maps $(E, \Phi)$ to the eigenvalues of $\Phi$ encoded in the characteristic polynomial $\det(\Phi-\eta \operatorname{Id})$ of $\Phi$, and is obtained as follows. Let $\{p_{1}, \ldots, p_k\}$ be a homogeneous basis for the algebra of invariant polynomials on the Lie algebra $\mathfrak{g}_{c}$ of $ \GC$, and let $d_{i}$ denote the degree of $p_i$. The {\it Hitchin fibration} is then given by
\begin{align} \label{hitchin_map}
\mathrm{Hit}\colon \ \mathcal{M}_{ \GC}&\longrightarrow \mathcal{A}_{ \GC}:=\bigoplus_{i=1}^{k}H^{0}\big(\Sigma,K^{d_{i}}\big),\\
(E,\Phi)&\mapsto (p_{1}(\Phi), \ldots, p_{k}(\Phi)),\nonumber
\end{align}
 where $\mathrm{Hit}$ is referred to as the {\it Hitchin~map}. It is a proper map for any choice of basis\footnote{In particular, it can be expressed in terms of the coefficients of $\det(\Phi-\eta \operatorname{Id})$.}, its generic fibers are abelian varieties, and makes the moduli space into a complex integrable system~\cite{hitchin87b}.

Each connected component of a generic fiber of the Hitchin map is an abelian variety. In the case of $\GC$-Higgs bundles this can be seen using spectral data \cite{BNR, hitchin87b}. Through the characteristic polynomial of the Higgs field of a $\GC$-Higgs bundle $(E,\Phi)$, one may define an algebraic curve, called the {\em spectral curve} of $(E,\Phi)$, which is generically smooth\footnote{When considering classical groups $\GC$, only for ${\rm SO}(2p,\C)$ one needs to consider a normalization of the curve.}:
 \begin{gather}S=\{{\rm det}(\Phi - \eta \operatorname{Id})=0\} \subset {\rm Tot}(K),\label{curva} \end{gather}
where ${\rm Tot}(K)$ is the total space of $K$, the map $\eta$ is the tautological section of $K$ on ${\rm Tot}(K)$, and by abuse of notation, we consider $\Phi$ as its pull-back to ${\rm Tot}(K)$ (the reader should refer to~\cite{BNR} for thorough details on the construction). We say that $(E,\Phi)$ lies in the {\em regular locus} of $\mathcal{M}_{\GC}$ if the curve~$S$ is non-singular, and denote the regular locus of the moduli space by $\mathcal{M}_{\GC}'$. Let $\pi \colon S \to \Sigma$ be the natural projection to $\Sigma$, and let $\eta \in H^0( S , \pi^*(K) )$ denote the restriction of the tautological section of $K$ to $S$. If $(E,\Phi)$ is in the regular locus, then there exists a line bundle $L \to S$ for which $E = \pi_* L$, and $\Phi$ is obtained by pushing down the map $\eta\colon L \to L \otimes \pi^*(K)$. In this way, one recovers the pair $(E,\Phi)$ from the pair $(S,L)$, which is referred to as the {\em spectral data} associated to the pair $(E , \Phi)$.

Note that the spectral curve $S$ of the pair $(E , \Phi)$ depends only on the characteristic polynomial of $\Phi$ and hence it only depends on the image of $(E,\Phi)$ under the Hitchin map. Therefore any point $a \in \mathcal{A}_{\GC}$ has an associated spectral curve $S_{a}$. If $a$ is in the regular locus of $\mathcal{A}_{\GC}$, in other words, if the associated spectral curve $S_{a}$ is smooth, then the spectral data construction identifies the fiber $\mathrm{Hit}^{-1}(a)$ of the Hitchin system with some subspace of $\mathrm{Pic}(S_{a})$, the Picard variety of the spectral curve $S$. The connected components of~$\mathrm{Pic}(S_{a})$ are isomorphic to copies of~$\mathrm{Jac}(S_{a})$, the Jacobian of $S_{a}$, and are labeled by the degree of the vector bundle~$E$ of a Higgs pair~$(E,\Phi)$. In particular, for $G_\C={\rm GL}(n,\C)$, the generic fibers are isomorphic to~$\mathrm{Jac}(S_{a})$, and one can see from here that the components of the regular fibers are abelian varieties. While much is known about the generic fibers of the $G_\C$-Hitchin fibration, there are still several interesting open questions. In particular, it would be interesting to understand the geometry of the generic fibers of the Hitchin fibration stated in the open problems of \cite{singapore}. For instance, it is interesting to consider the following:

 \begin{open} Considering the notion of ``strong real form''\footnote{The notion of strong real form is a refinement of the notion of real form. For example for ${\rm SL}(2,\C)$, there are three equivalence classes of strong real forms corresponding to ${\rm SU}(2, 0),~ {\rm SU}(0, 2)$
and ${\rm SU}(1, 1)$.} from {\rm \cite{adam}}, describe the corresponding Higgs bundles and determine which ones define singular spectral curves. \end{open}

When considering arbitrary groups $\mathcal{G}$, the algebraic curve defined by the characteristic polynomial of $\Phi$ is not always generically smooth (for real groups, see for instance the case of $\mathcal{G}={\rm SU}(p,q)$ and the spectral data described in~\cite{thesis}). In this case, one may consider cameral covers \cite{Don95}: these are $K$ valued covers of $\Sigma$ with an action of~$\mathcal{W}$, the Weyl group of $\mathcal{G}$. The fiber of the associated fibration can be described in terms of these covers, and over a generic point of the base the cover is a $\mathcal{W}$-Galois cover (e.g., see \cite{Don93, Don95, DG02, Fal93, Sco98}, and \cite{donagipantev} for further references). In this set up, there is a natural discriminant locus in the Hitchin base, away from which the connected component of the fiber is isomorphic to a certain abelian variety which can be described as a generalized Prym variety of the cameral cover. However, the study of the singular fibers, even from the perspective of cameral covers, is not fully understood.

 \begin{open} Give a comparison of what is known for singular fibers of the Hitchin fibration from the perspective of cameral covers and of spectral data.
 \end{open}

Cameral covers have shown to be very useful tools to understand the moduli spaces of principal Higgs bundles and their relation to many other fields. However, the abstraction of the method and the constructions of the covers can sometimes make certain properties of the moduli spaces very difficult to discern. Although most objects are defined in the above papers in a general way, their description and study for particular groups is still being done by many researchers (e.g., see recent developments for real Higgs bundles in~\cite{ana2,ana1}, where interesting comparisons with classical spectral data are carefully explained).

\begin{open}Extend the cameral cover methods of {\rm \cite{ana2}} for ${\rm SU}(p,p+1)$-Higgs bundles to all other real Higgs bundles which lie completely over the singular locus of the Hitchin fibration, i.e., to those Higgs bundles whose characteristic polynomial defines singular curves through~\eqref{curva}.
\end{open}

\subsection{The singular locus of the Hitchin fibration}
As mentioned before, the fiber $\mathrm{Hit}^{-1}(a)$ over $a \in \mathcal{A}_{\GC}$ is said to be singular when the corresponding spectral curve $S_a$ defined as in~\eqref{curva} is singular. The most singular fiber is the nilpotent cone\footnote{The name was given by Laumon \cite{laumon}, to emphasize the analogy with the nilpotent cone in Lie algebra.}, which sits over $\mathbf{0} \in \mathcal{A}_{\GC}$. One of the tools to study the nilpotent cone is the moment map~$\mu$ of the $S^1$ action $(E, \Phi, h) \to (E, \mathrm{e}^{i\theta} \Phi, h)$ \cite{hauseldiss, hitchin87}. Moreover, the nilpotent cone is preserved by the flow by~$\mu$ \cite[Theorem~5.2]{Hausel:1998aa}, and since points of $\M_{\GC}$ flow towards the nilpotent cone, it encodes the topology of the moduli space.

The nilpotent cone has primarily been studied for ${\rm SL}(n,\C)$ and ${\rm GL}(n,\C)$, and much of its geometry remains unknown for the moduli spaces of $G_\C$-Higgs bundles. For ${\rm SL}(n,\C)$ and ${\rm GL}(n,\C)$-Higgs bundles, the irreducible components of the nilpotent cone are labeled by connected components of the fixed point set of the $S^1$ action. Among these components is the moduli space $\mathcal{N}$ of semistable bundles\footnote{Given a stable bundle $E$, take $\Phi=0$, then the Higgs bundle $(E, 0)$ is stable and trivially fixed by the $S^1$-action.}.

Other singular fibers have been the subject of more recent research (e.g., see \cite{orthogonal1,gothenoliveira,nonabelian}). In particular, in the case of ${\rm GL}(n,\C)$-Higgs bundles while when the spectral curve $S$ is smooth, the corresponding fiber $\mathrm{Hit}^{-1}(S)$ can be identified with the Jacobian $\mathrm{Jac}(S)$ of all line bundles $L \rightarrow S$ of degree $0$, when the spectral curve $S$ is not smooth, the corresponding fiber is seen to be the \emph{compactified Jacobian} \cite{BNR, schaub} (see also \cite[Fact 10.3]{melo} for a clear explanation). The compactified Jacobian $\overline{\mathrm{Jac}}(S)$ is the moduli space of all torsion-free rank-1 sheaves on $S$, where the usual ``locally-free'' condition is missing. Moreover, when $S$ is not integral, the fine moduli space needs to be considered. A more intuitive definition of the compactified Jacobian is the following: consider a path of smooth curves $S_t$ approaching a singular curve $S_0$; since the limit of $\mathrm{Jac}(S_t)$ does not depend on the choice of smooth family \cite{igusa}, this limit is $\overline{\mathrm{Jac}}(S)$.

In the case of ${\rm SL}(2,\C)$-Higgs bundles, much work has been done on the singular fibers of the corresponding Hitchin fibration. For example, see \cite[Section~5.2.2]{frenkelwitten} for connectedness of the fiber of $\M_{{\rm SL}(2,\C)}$ when $S$ is irreducible and has only simple nodes; see \cite{gothenoliveira} for a fuller description of the singular fibers\footnote{The authors of \cite{gothenoliveira} actually study the slightly more general situation of ``twisted Higgs bundles'' where the canonical bundle $K$ is replaced by a line bundle $L$ with $\deg (L)>0$. For this setting of $L$-twisted ${\rm SL}(2,\C)$-Higgs bundles, the monodromy action was considered in \cite{monodromy1}.} of $\M_{{\rm SL}(2,\C)}$, and see \cite{monodromy2} for the monodromy action around singular fibers of the Hitchin fibration.

\begin{open} Building on the results for ${\rm SL}(2,\C)$-Higgs bundles, describe the singular fibers of the Hitchin fibration for arbitrary $G_\C$.
\end{open}

When the singular fiber lies above some particular types of spectral curves, one may describe the fibers by considering some modified version of spectral data, leading to the following natural question:

\begin{open}Extending on {\rm \cite{nonabelian}}, obtain a geometric description of the fibers of the Hitchin fibration of $G_\C$-Higgs bundles which lie over points of the Hitchin base defining curves through \eqref{curva} with equation $\det(\Phi-\eta \operatorname{Id})=P^k(\eta)$ for $k\geq2$, and for which $\{P(\eta)=0\}$ is generically smooth $($and thus defines itself a smooth spectral curve$)$.
\end{open}

When the spectral curve has defining equation $\det(\Phi-\eta \operatorname{Id})=P^2(\eta)$, components of the fiber were studied in \cite{nonabelian}, and the full fibers of the Hitchin fibration with that base point are described in~\cite{lucas}. From a different perspective, in terms of fiber products of spectral curves, certain singular spectral curves were considered in~\cite{isogenies}. While not much is known about the singular fibers of the Hitchin fibration for $G_{\C}$-Higgs bundles, one may deduce properties of the whole moduli space by considering the monodromy action of the natural Gauss--Manin connection of the fibration. In the case of ${\rm SL}(2,\C)$-Higgs bundles, the study of the monodromy was done in~\cite{monodromy2}, where an explicit formula was used to understand connectivity of the moduli space. The work was later extended to twisted rank~2 Higgs bundles in \cite{monodromy1}, and to all ${\rm SL}(n,\C)$-Higgs bundles in \cite{monodromy3}. However, the general understanding of the monodromy action for other groups remains open:
\begin{open}
Give a geometric description of the monodromy action for the Hitchin fibration of $G_\C$-Higgs bundles.
\end{open}

Finally, since the moduli spaces $\M_{\GC}$ are often not smooth, it is important to understand the singularities of $\M_{\GC}$. For a beautiful survey on recent developments in the theory of moduli spaces of sheaves on projective varieties, and implications for Higgs bundles, the reader may refer to~\cite{Carlos16}. In the case of parabolic Higgs bundles, a description of the Hitchin fibration was given recently in~\cite{david17}, and it would be very interesting to understand the above considerations and open questions in this other setting.

\section{Higgs bundles and limiting structures} \label{sec-limmits}

Many conjectures from mathematics and physics about $\M_{\GC}$ remain open because they require a finer knowledge of the ends of the moduli space than what is provided by traditional algebro-geometric techniques. In this section we shall restrict our attention to ${\rm SL}(2,\C)$-Higgs bundles, and note that for other groups most of the questions mentioned here remain open. In this setting, one has the following conjecture of Hausel:

\begin{Conjecture}[{\cite[Conjecture 1]{hauseldiss}}] There are no non-trivial $L^2$ harmonic forms on the Hitchin moduli space.
\end{Conjecture}

There is a similar conjecture for the moduli space of monopoles which is called the Sen Conjecture. By analogy, the conjecture for the Hitchin moduli space is sometimes called the Sen Conjecture as well. In order to obtain a finer knowledge of the ends of $\M_{{\rm SL}(2,\C)}$, finer descriptions of solutions of Hitchin's equations near the ends are needed. A number of recent results~\cite{MSWW15,MSWW14, MSWW17, Mochizuki:2015aa} demonstrate the power of constructive analytic techniques for describing the ends of the Hitchin moduli space.

Fixing a stable Higgs bundle $(E, \Phi)$ in $\M_{{\rm SL}(2,\C)}$, the ray of Higgs bundles with harmonic metric $(E, t \Phi, h_t)$ approaches the ends of the moduli space as $t \to \infty$. In order to understand what the behavior of the harmonic metrics $h_t$ as $t \to\infty$ is, note that in the limit the curvature $F_{D(\delbar_E, h_t)}$ concentrates at the ramification points $Z \subset \Sigma$ of $\pi\colon S \rightarrow \Sigma$ and vanishes everywhere else. The decay is exponential in $t$, leading to the following result:
\begin{Theorem}[{\cite[Theorem 2.7]{Mochizuki:2015aa}}] On a compact subset $\overline{U}$ of $\Sigma-Z$, there exist
positive constants $c_0$ and $\eps_0$ such at any point in $\overline{U}$
\begin{gather*} \big|[\varphi, \varphi^{\dagger_{h_t}} ] \big|_{h_t, g_\Sigma} \leq c_0 \exp(-\eps_0 t). \end{gather*}
\end{Theorem}

Consequently, the limiting hermitian metric is singular at the ramification points $Z \subset \Sigma$ and
\begin{gather} \label{eq:hitchindecoupled}
 F_{D(\delbar_E, h_\infty)}=0, \qquad [\Phi, \Phi^{*h_\infty}]=0 , \qquad \delbar_E \Phi =0.
\end{gather}
It is often said that Hitchin's equations ``abelianize'' asymptotically. The vanishing of the Lie bracket in~\eqref{eq:hitchindecoupled} reflects the deeper expectation that the metric $h_\infty$ is the pushforward of a~singular harmonic metric~$h_L$ on the spectral data $L \rightarrow S$. This has been proved when $S$ is smooth, i.e., $(E, \Phi) \in \M'_{{\rm SL}(n,C)}$, by Mazzeo--Swoboda--Weiss--Witt when $n=2$ \cite{MSWW14} and generalized to any rank by Fredrickson~\cite{FredricksonSLn}. In \cite[Theorem~5.1]{Mochizuki:2015aa} Mochizuki proves this for all of $\mathcal{M}_{{\rm SL}(2,\C)}$, making no assumptions about the smoothness of~$S$.

\subsection{The ends of the regular locus} More is known about the ends of the regular locus $\M'_{{\rm SL}(2,\C)}$. For $t$ large but finite, the harmonic metric $h_t$ is close to an approximate harmonic metric $h_t^{\mathrm{approx}}$, constructed by desingulari\-zing~$h_\infty$~\cite{MSWW14}. On small disks around points in $Z \subset \Sigma$, the approximate metric is equal to a~smooth local model solution
\begin{gather*}
 \delbar_E=\delbar, \qquad t\Phi = t\begin{pmatrix}0 & 1 \\ z & 0 \end{pmatrix} \de z, \qquad
 h_t^{\mathrm{model}} = \begin{pmatrix} |z|^{1/2} \e^{u_t(|z|)} & \\ & |z|^{-1/2} \e^{-u_t(|z|)} \end{pmatrix},
\end{gather*}
where $u_t(|z|)$ comes as a solution of a $t$-rescaled Painlev\'e III ODE with boundary conditions given by $u_t(|z|) \sim -\frac{1}{2} \log(|z|)$ near $|z|=0$ (so that $h_t$ is smooth), and $\lim\limits_{|z| \to \infty} u_t(|z|)=0$. Note that in this same local gauge, the singular limiting metric $h_\infty$ would be equal to
\begin{gather*}
 h_\infty = \begin{pmatrix} |z|^{1/2} & \\ & |z|^{-1/2} \end{pmatrix}.
\end{gather*}
Outside of small disks around $Z \subset \Sigma$, the approximate harmonic metric $h_t^{\mathrm{approx}}$ is equal to $h_\infty$.

The approximate description of $h_t$ by $h_t^{\mathrm{approx}}$ has already been useful in \cite{MSWW17} for describing the hyperk\"ahler metric on $\M'_{{\rm SL}(2,\C)}$ near the ends. There are two natural hyperk\"ahler metrics on $\M'_{{\rm SL}(2,\C)}$: first, the hyperk\"ahler metric $g_{L^2}$ on $\M_{{\rm SL}(2,\C)}$ restricts to $\M'_{{\rm SL}(2,\C)}$; second, there is a metric $g_{\mathrm{sf}}$ on $\M'_{{\rm SL}(2,\C)}$, known as the semiflat metric because $g_{\mathrm{sf}}$ is flat on the half-dimensional torus fibers~\cite{freedsemiflat}. The metric $g_{L^2}$ comes from taking the $L^2$ metric on triples $(E, t \Phi, h_t)$, while the semiflat metric $g_{\mathrm{sf}}$ comes from taking the $L^2$ metric on the moduli space of triples~$(E, t \Phi, h_\infty)$~\cite{MSWW17}. Consequently, Mazzeo--Swoboda--Weiss--Witt are able to describe the difference between~$g_{L^2}$ and~$g_{\mathrm{sf}}$ using their careful description of~$h_t$ and $h_\infty$. They prove
 \begin{Theorem}[{\cite[Theorem 1.2]{MSWW17}}]
 The metric $g_{L^2}$ admits an asymptotic expansion
 \begin{gather} \label{eq:MSWWmetric}
 g_{L^2} = g_{\mathrm{sf}} + \sum_{j=0}^{\infty} t^{\frac{4-j}{3}}G_j + \mathcal{O}\big(e^{-\beta t}\big)
 \end{gather}
as $t \to \infty$. Here each $G_j$ is a symmetric two-tensor.
 \end{Theorem}

 \begin{open}\label{pool} Are the polynomial correction terms $G_j$ non-zero?
 \end{open}

Separate from this description of the ends with PDE techniques, a remarkable conjectural picture of the asymptotic geometry of $\M_{{\rm SL}(n,\C)}$ has emerged from physics in the work of Gaiotto, Moore, and Neitzke \cite{GMNwallcrossing, Neitzkehyperkahler}. Their starting point is the semiflat metric $g_{\mathrm{sf}}$ on $\M'_{{\rm SL}(n,\C)}$ which is too homogeneous to extend to all of $\M_{{\rm SL}(n,\C)}$. They give a recipe for constructing a complete hyperk\"ahler metric $g_{\mathrm{GMN}}$ on $\M_{{\rm SL}(n,\C)}$ differing from $g_{\mathrm{sf}}$ by ``quantum corrections'' which are computed by counting certain BPS states in supersymmetric field theory. In particular, the quantum corrections have the following size
 \begin{gather}\label{eq:GMN}
 g_{\mathrm{GMN}} = g_{\mathrm{sf}} + O \left( \sum_{\gamma \in H_1(S_a, \Z)} \Omega(\gamma, a) \mathrm{e}^{-t \big| \int_\gamma \eta\big|} \right)
 \end{gather}
 as $t \to \infty$. In this formula, as in the previous sections, $a$ is a point in $\cA_{{\rm SL}(n,\C)}$, the corresponding spectral
 curve is $S_a$ with tautological one-form $\eta$; and the sum is over all loops $\gamma$ in $S_a$. The~$\Omega(\gamma, a)$ are BPS counts in supersymmetric field theory. These are $\Z$-valued and piecewise-constant, jumping across certain walls in the parameter space. The jumps are constrained to satisfy the Kontesevich--Soibelman wall-crossing formula~\cite{kontsevichsoibelman}, and thus $g_{\mathrm{GMN}}$ is smooth.
 Moreover, Gaiotto--Moore--Neitzke conjecture
 \begin{Conjecture}[\cite{GMNwallcrossing}]\label{GMN}
 The hyperk\"ahler metric $g_{\mathrm{GMN}}$ on the moduli space is the natural hyperk\"ahler metric $g_{L^2}$ on $\M_{{\rm SL}(n,\C)}$.
 \end{Conjecture}
If the conjecture of Gaiotto--Moore--Neitzke is correct, then all of the symmetric two-ten\-sors~$G_j$ appearing in~\eqref{eq:MSWWmetric} vanish, answering Open Question~\ref{pool}. Note that there is already evidence that this happens on the Hitchin section~\cite{DumasNeitzkeprep}. As $a \in \mathcal{A}_{{\rm SL}(n,\C)}$ approaches the singular locus $\cA^{\mathrm{sing}}_{{\rm SL}(n,\C)}$, the spectral curves~$S_a$ become singular. In particular, there is at least one loop~$\gamma_0$ on~$S_a$ which pinches; hence the ``quantum correction'' in~\eqref{eq:GMN} corresponding to $\gamma_0$ is not exponentially suppressed as $t$ approaches~$\infty$. While we have focused here on ${\rm SL}(2,\C)$-Higgs bundles, equivalent questions (and conjectures) may be asked for more general $G_\C$-Higgs bundles, providing several new lines of research:

 \begin{open}Generalize the above results to the case of $G_\C$-Higgs bundles for arbitrary~$G_\C$.
 \end{open}

\section{Higgs bundles on singular curves} \label{sec-curves}

While we have considered before $G_\C$-Higgs bundles on a compact Riemann surface $\Sigma$, principal Higgs bundles can also be defined over singular spaces $X$, and in particular, over singular curves. For simplicity, we shall begin by considering nodal curves $X$, i.e., irreducible projective curves~$X$ whose singularities are nodes. In order to generalize the notion of vector bundles on smooth projective curves, one may consider the torsion-free sheaves on the nodal curve~$X$. Through the work of Bhosle in~\cite{bhosle1}, the category of torsion-free sheaves on a nodal curve $X$ and category of generalized parabolic bundles over its normalization $\widetilde X$ are equivalent. Moreover, a first general construction of compactified moduli spaces for semistable $G_\C$-bundles on an irreducible complex projective curve $X$ with exactly one node was given in~\cite{alexander}.

\subsection[Singular principal $G_{\C}$-Higgs bundles]{Singular principal $\boldsymbol{G_{\C}}$-Higgs bundles} Through the work of \cite{alexander}, it was shown in \cite{andrea} that one can treat a principal $G_\C$-Higgs bundle over a nodal curve $X$ as a particular type of vector bundle on the normalization of the curve called a {\it descending bundle}, objects which are in one-to-one correspondence with the following {\it singular principal $G_\C$-Higgs bundles}.

 \begin{Definition}
 A singular principal $G_\C$-Higgs bundle is a triple $(\mathcal{E},\tau, \Phi)$ where
 \begin{itemize}\itemsep=0pt
 \item $\mathcal{E}$ is a locally free sheaf;
 \item $\tau\colon {\rm Sym}^*(\mathcal{E}\otimes V)^{G_\C}\rightarrow \CO_X$, for a fixed faithful representation $\GC \rightarrow {\rm GL}(V)$ of $\GC$;
 \item $\Phi\colon X\rightarrow {\rm End}(\mathcal{E})\otimes \Omega_X^1$ is a section;
\end{itemize}
 \end{Definition}

\begin{open} Give, if possible, a notion of the Hitchin fibration for these singular principal $G_\C$-Higgs bundles, and describe the geometry of the smooth and singular fibers.
 \end{open}

 While Schmitt \cite{alexander} and Bhosle \cite{bhosle1} proved that there is a moduli space of singular $G_\C$-bundles on $X$ with good specialization properties, Seshadri gave a further study of the spaces of torsion-free sheaves on nodal curves and generalizations to, among others, ramified $G_\C$-bundles~\cite{Seshadri}. When considering degenerations of moduli spaces of vector bundles on curves, which are closely related to the singular $G_\C$-bundles mentioned above, the reader may want to consider the conjectures presented in \cite{seshadri2}.

Just as one may define parabolic Higgs bundles on Riemann surfaces to consider Higgs bundles on marked curves, one may extend these objects to singular curves. Extending the notion of a~parabolic vector bundle on a~smooth curve, Bhosle defined generalized parabolic sheaves (GPS) on any integral projective curve~$X$~\cite{bhosle2} and generalized parabolic bundles (GPB)~\cite{bhosle4}. She constructed the moduli spaces of GPS and GPB, and studied the correspondences appearing when curves $X$ are obtained from blowing up finitely many nodes in a space~$Y$. Moreover, Bhosle also extended the notion of parabolic Higgs bundles to that of generalized parabolic Higgs bundles (GPH) on the normalization~$X$ of an integral projective curve $Y$ \cite{bhosle2,bhosle4,bhosle3}. In particular, she constructed a birational morphism from the moduli space of {\it good} GPH on $X$ to the moduli space of Higgs bundles on~$Y$, and defined a proper Hitchin map on the space of GPH. In this context, the following question is natural.

\begin{open}Generalize {\rm \cite{bhosle3}} to define generalized parabolic $G_\C$-Higgs bundles on $X$, as well as a Hitchin fibration.
 \end{open}

Moreover, the open questions mentioned for the moduli spaces of classical $G_\C$-Higgs bundles may also be considered both for parabolic Higgs bundles on Riemann surfaces~$\Sigma$, and for parabolic Higgs bundles on integral projective curves~$X$. In particular, Bhosle studied recently the relationship between Higgs bundles and the compactified Jacobian of a spectral curve \cite{bhosle3}, considering Higgs bundles on the normalization $X$ of integral projective curves~$Y$, leading to an analogous question to that stated for classical Higgs bundles:
 \begin{open}
 Obtain a geometric description of the singular fibers of the Hitchin fibration for generalized parabolic $G_\C$-Higgs bundles on $X$.
 \end{open}

The study of Higgs bundles on singular curves may also be considered in a limiting setting, where one begins with Higgs bundles on a smooth curve and parametrically tunes the curve to degenerate to a singular curve. The particular case of vector bundles on smooth curves degenerating to an irreducible curve with one double point was considered in~\cite{ivan}. The case of the degeneration of the moduli space of Higgs bundles on smooth projective curves when the curve degenerates to an irreducible curve with a single node was studied in~\cite{Balaji}.

\begin{open} Obtain equivalent degenerations to those in {\rm \cite{Balaji,ivan}} for the moduli spaces of $\GC$-Higgs bundles.
 \end{open}

In particular, as explained in \cite{Balaji}, their degeneration is analogous to the models constructed by Gieseker and Nagaraj--Seshadri for the case of the moduli spaces for which the Higgs structure is trivial. It should be noted that in \cite{Balaji} the authors also construct a~corresponding canonical relative proper Hitchin map, whose fiber provides a new compactification of the Picard variety of smooth curves with normal crossing singularities. In their setting, the single node on the base curve leads to an {\it irreducible vine curve} with $n$-nodes appearing as the spectral curve. It would then seem natural that the quasi-abelianization of~\cite{Balaji} (resembling Hitchin's classical abelianization) could be generalized.
\begin{open}Describe the quasi-abelianization of the moduli space of $G_\C$-Higgs bundles for different degenerations, following the techniques of~{\rm \cite{Balaji}}.
\end{open}
Finally, since it is important to find natural compactifications of open moduli, and torsion-free sheaves on nodal curves play an important role in
 \cite{rahul} within the study of the compactification of the universal moduli space of slope-semistable vector bundles over the compactification $\overline{M}_g$ of the moduli space of genus $g$ curves, it is natural to ask the following:

\begin{open} Understand the relation between the degenerations of moduli spaces of Higgs bundles above, and the known compactifications of $\mathcal{M}_{G_\C}$.
\end{open}

 Since one has the correspondence between Langlands dual groups $G_\C$ and $^LG_\C$, once the corresponding moduli spaces are understood, and the Hitchin fibrations are shown to exist, one may also want to consider, if possible, the duality between the fibrations. In particular, the work of Arinkin~\cite{dima} for rank 2 Higgs bundles on the auto-duality of compactified Jacobians for curves with plane singularities would allow one to understand both this setting, as well as the one of singular fibers of the classical Hitchin fibration, leading to an intermediate question:
 \begin{open}
 Extend the constructions of {\rm \cite{dima}} to the setting of generalized parabolic $G_\C$-Higgs bundles \`a la Bhosle {\rm \cite{bhosle3}}.
 \end{open}

\section{Higgs bundles and branes within singular fibers} \label{sec-branes}

The appearance of Higgs bundles (and flat connections) within string theory and the geometric Langlands program has led researchers to study the {\it derived category of coherent sheaves} and the {\it Fukaya category} of these moduli spaces. Therefore, it has become fundamental to understand Lagrangian submanifolds of the moduli space of Higgs bundles supporting holomorphic sheaves ($A$-branes), and their dual objects ($B$-branes). For $^L G_\C$ the Langlands dual group of $G_\C$, there is a correspondence between invariant polynomials for $G_\C$ and $^L G_\C$ giving an identification $\mathcal{A}_{G_\C} \simeq \mathcal{A}_{^L G_\C}$ of the Hitchin bases.

\subsection{Construction of branes}

Through the Hitchin fibrations, the two moduli spaces $\mathcal{M}_{G_\C}$ and $\mathcal{M}_{^L G_\C}$ are then torus fibrations over a common base and their non-singular fibers are dual abelian varieties \cite{donagipantev, hauselthaddeus}, answering some of the conjectures presented in~\cite{SYZ}. Kapustin and Witten give a physical interpretation of this in terms of S-duality, using it as the basis for their approach to the geometric Langlands program \cite{kapustinwitten}. In this approach a crucial role is played by the various types of branes and their transformation under mirror symmetry. Adopting the language of physicists, a Lagrangian submanifold of a symplectic manifold supporting a hyperholomorphic sheaf is called (the base of) an {\em A-brane}, and a complex submanifold supporting a hyperholomorphic sheaf is (the base of) a {\em B-brane}. A~submanifold of a hyperk\"ahler manifold may be of type $A$ or $B$ with respect to each of the complex or symplectic structures, and thus choosing a triple of structures one may speak of branes of type $(B,B,B)$, $(B,A,A)$, $(A,B,A)$ and $(A,A,B)$\footnote{One should note that since the complex structures satisfy the quaternionic equations, and the symplectic forms are obtained through them, branes of types $(A,A,A)$, $(A,B,B)$, $(B,A,B)$, $(B,B,A)$ do not exist.}. Throughout these notes we shall follow the convention in~\cite{kapustinwitten} and fix the three complex structures $I$, $J$ and $K$, such that $I$ is induced from the Riemann surface~$\Sigma$, and~$J$ from the complex group~$\GC$.

It is hence natural to seek constructions of different families of branes inside the moduli space $\mathcal{M}_{G_c}$, understand their appearance within the Hitchin fibration, and describe their mirror families of branes. In the context of Higgs bundles, branes were first considered by Kapustin and Witten in 2006 in \cite{kapustinwitten}, where much attention was given to the $(B,A,A)$-brane of $G$-Higgs bundles inside $\M_{G_{\mathbb{C}}}$, where $G$ is a real form of the complex Lie group~$G_\C$. Soon after, examples of brane dualities were considered in \cite{sergei1}; in particular the case of $G$-Higgs bundles for compact real forms $G$ of low rank was considered, in which case the $(B,A,A)$-brane lies completely inside the nilpotent cone. While partial results exist for these branes over singular fibers, the more global picture remains unknown.

 \begin{open} Give a geometric description of all $(B,A,A)$-branes of $G$-Higgs bundles which live completely inside the most singular fiber of the Hitchin fibration, the nilpotent cone of~$\mathcal{M}_{G_\mathbb{C}}$.
 \end{open}

 The study of branes within moduli spaces of Higgs bundles continued evolving slowly, until natural generic methods to construct families of all types of branes in $\M_{\GC}$ were introduced in~\cite{slices}. These branes were constructed as fixed point sets of certain families of involutions on the moduli spaces of complex Higgs bundles. Consider $\sigma$ an anti-holomorphic involution fixing a real form $G$ of $\GC$, and $\rho$ the anti-holomorphic involution fixing the compact real form of $\GC$. Then, through the Cartan involution $\theta=\sigma \circ \rho$ of a real form~$G$ of~$G_\C$, one may define
\begin{gather*}
i_1\big(\bar \partial_E, \Phi\big):=\big(\theta(\bar \partial_E\big),-\theta( \Phi)).
\end{gather*}
 Moreover, a real structure $f \colon \Sigma \to \Sigma$ on $\Sigma$ induces an involution
\begin{gather*}
i_2(\bar \partial_E, \Phi):=\big(f^*( \partial_E),f^*( \Phi^* )\big)= \big(f^*\big(\rho\big(\bar \partial_E\big)\big), -f^*( \rho(\Phi) )\big).
\end{gather*}
 Lastly, by setting $i_3 := i_1 \circ i_2$, one may define a third involution:
\begin{gather*}
i_3(\bar \partial_E, \Phi)=\big(f^* \sigma\big(\bar \partial_E\big),f^*\sigma( \Phi)\big).
\end{gather*}
The fixed point sets of the induced involutions $i_1$, $i_2$, $i_3$ introduced in \cite{slices} are branes of type $(B,A,A),$ $(A,B,A)$ and $(A,A,B)$ respectively, and through the associated {\it spectral data} their topological invariants can be described using $KO$, $KR$ and equivariant $K$-theory. In particular, it was shown that among the fixed points of $i_1$ are solutions to the Hitchin equations with holonomy in $G$. Moreover, those fixed by $i_2$ were shown to give real integrable systems, fibered as a Lagrangian fibration over a real slice of the Hitchin base \cite{aba}. In order to construct the fourth type of branes, $(B,B,B)$-branes, one may consider Higgs bundles for a complex subgroup of $G_{\C}$, but these branes would not appear through a symmetry in the spirit of the above constructions\footnote{An instance of this setting was recently explored in \cite{emilio4}, where Langlands duality was studied for branes appearing though Borel subgroups.}. On the other hand, it is shown in \cite{cmc} that one may construct $(B,B,B)$-branes by considering the subspaces of $\Gamma$-equivariant Higgs bundles for $\Gamma$ a finite group acting on the Riemann surface~$\Sigma$. In particular, it was shown in~\cite{cmc} that for $G_\C={\rm SL}(2,\C)$, these branes would be mid-dimensional only under very restrictive conditions, and no equivalent result has been shown for higher rank groups.

 \begin{open} Describe the mid-dimensional $(B,B,B)$-branes appearing through $\Gamma$-equivariant Higgs bundles when $\Gamma$ is a group of any rank, and classify those components completely contained in the singular locus of the Hitchin fibration.
 \end{open}

 The construction of branes following the procedures of \cite{aba,slices} have recently been generalized to the space of framed instantons \cite{emilio1}, Higgs bundles over K3 surfaces \cite{emilio2}, Higgs bundles over elliptic curves \cite{emilio3}, quiver varieties \cite{vic,vic2}, more general hyperk\"ahler spaces \cite{indranil11}, and principal Schottky bundles \cite{carlos}. Moreover, many of the geometric properties of the branes in \cite{aba, slices} are yet unknown, and researchers continue to study them (e.g., see \cite{tom, indranil11}). In the case of finite group actions, the branes introduced in \cite{cmc} were later studied in \cite{copy1} from the perspective of character varieties, and many of their properties remain unknown.

 \begin{open}In the spirit of {\rm \cite{biswas1}} and {\rm \cite{biswas2}}, describe the Brauer groups and automorphism groups of the branes mentioned above.
 \end{open}

 Finally, it should be mentioned that in the last couple of years researchers have found other novel ways in which branes can be constructed within the moduli space of Higgs bundles, and which are yet to be generalized to other settings. Examples of these are Nahm branes~\cite{emilio5}, branes appearing through spinors~\cite{hitchin1}, through moment maps~\cite{davide1}, and through Borel subgroups \cite{emilio4}. However, since Lagrangian branes can appear in any of three types, it is of interest to understand families of Lagrangian branes of each types which are related in some geometric fashion, but it is not yet known of other triples of families of branes appearing within hyperk\"ahler spaces other than the ones obtained thought the methods of \cite{slices}.

 \begin{open}Construct natural triples of families of branes in \MGC, and more ge\-ne\-rally, within hyperk\"ahler spaces.
 \end{open}

Considering the appearance of branes through real structures on Riemann surfaces, one should also be able to impose other structures on the surfaces to construct novel branes. A canonical example of such structure would be that of a {\it log-symplectic structure} on $\Sigma$, also called a {\it b-Poisson structure}. These structures are given by Poisson structures $\pi\in \mathfrak{X}^{2}(\Sigma)$ for which $\pi$ has only non-degenerate zeros. In particular, $\pi$ is generically symplectic.
These structures were completely classified by O. Radko \cite{radko}, where she noted that every surface (orientable or not) has a log-symplectic structure. The sets of invariants of log-symplectic structures are:
\begin{itemize}\itemsep=0pt
 \item The zero curves $\gamma_1, \ldots, \gamma_n$, taken with orientation defined by $\pi$;
\item The periods associated to each $\gamma_i$;
\item The volume invariant of $\pi$.
\end{itemize}Since there is a natural relation between the data defining log-symplectic structures $\pi$ on $\Sigma$ as in~\cite{radko}, and the real structures $f\colon \Sigma\rightarrow \Sigma$ considered in~\cite{aba} to define $(A,B,A)$-branes, a~natural question is the following.

\begin{open}Which branes of Higgs bundles are characterized by log-symplectic structures, and how do these relate to the $(A,B,A)$-branes introduced in~{\rm \cite{aba,slices}}?
\end{open}

\subsection{Langlands duality} While it is understood that Langlands duality exchanges brane types, the exact correspondence is not yet known\footnote{Strictly speaking, Langlands is a correspondence between local systems on $\Sigma$ (or more precisely, coherent sheaves on the moduli space), and $D$-modules over the moduli stack Bun of bundles on $\Sigma$.}. As mentioned before, the first instances of the correspondence being studied for low rank Higgs bundles appeared in \cite{kapustinwitten} and \cite{sergei1}, but no proof has yet been given of a~pair of branes of Higgs bundles being dual. In the case of $(B,A,A)$-branes of $G$-Higgs bundles, it was conjectured in \cite{slices} how the duality should appear:

 \begin{Conjecture}[{\cite[Section 7]{slices}}] \label{duality} 
 The support of the dual $(B,B,B)$-brane in $\mathcal{M}_{^LG_\C}$ to the $(B,A,A)$-brane $\MG\subset \mathcal{M}_{G_{\C}}$ is the moduli space $\mathcal{M}_{\check{H}}\subset \mathcal{M}_{^L \GC}$ of $\check{H}$-Higgs bundles for $\check{H}$ the group associated to the Lie algebra $\check{\mathfrak{h}}$ in {\rm \cite[Table~1]{nadler}}.
 \end{Conjecture}

 Support for this conjecture is given in \cite{clases} for the group $G=U(m,m)$ by considering the spectral data description of the brane in \cite{umm}, and in \cite{orthogonal1,orthogonal2} for the groups $G={\rm SO}(p+q,p)$ and ${\rm Sp}(2p+2q,2p)$. One should note that, in contrast with the $(A,B,A)$ and $(A,A,B)$ branes considered in~\cite{slices}, for any $q>1$ the $(B,A,A)$-branes studied in \cite{orthogonal1} lie completely over the singular locus of the Hitchin fibrations. For these branes of orthogonal Higgs bundles, support for Conjecture~\ref{duality} is obtained from the description of how the brane intersects the most generic fibers of the Hitchin fibration: indeed the rank of the hyperholomorphic sheaf depends on the number of components in this intersection, which remains constant for different $q$, leading to the following conjecture.

 \begin{Conjecture}[{\cite[Section 8]{orthogonal1}}] For all $q$ even $($and for all $q$ odd$)$, the $(B,A,A)$-brane of ${\rm SO}(p+q,p)$-Higgs bundles has dual $(B,B,B)$-brane obtained by considering the same base of Conjecture~{\rm \ref{duality}} and the same hyperholomorphic bundle supported on it\footnote{This hyperholomorphic bundle being, for example, the one introduced by Hitchin in~\cite{clases}.}, and it is only the way in which these spaces are embedded into the different Langlands dual moduli spaces which depends on~$q$.
 \end{Conjecture}

In particular, the support of branes for $q$ odd and $q$ even are dual to each other as hyperk\"ahler moduli spaces of complex Higgs bundles.
From the description of the invariant polynomials appearing for $(B,A,A)$-branes of $G$-Higgs bundles in \cite{thesis}, one can see that the majority of these branes lie over the singular locus of the Hitchin fibration. However, it is still possible to describe the intersection of these branes with the most regular of singular fibers. For example, it was shown in \cite{nonabelian} that the generic intersections of
 the $(B,A,A)$-branes of ${\rm SL}(m,\mathbb{H})$, ${\rm SO}(2n,\mathbb{H})$ and ${\rm Sp}(2m,2m)$-Higgs bundles with the fibers of the Hitchin fibrations are not abelian varieties, but are instead moduli spaces of rank 2 bundles on a spectral curve, satisfying certain natural stability conditions. In order to fully understand Langlands duality for branes one would need to understand how different branes in a moduli space intersect, and thus the particular case of branes within the nilpotent cone is of much importance.

\begin{open}
Describe the intersections and relations between all $(B,A,A)$-branes in the nilpotent cone of $G_\C$-Higgs bundles, for $G_\C$ an arbitrary group.
\end{open}

While a first step towards an answer would be to consider $(B,A,A)$-branes of $G$-Higgs bundles, or those constructed in the papers mentioned above, a more general perspective considering generators of the corresponding Fukaya category would be ideal. For a short review of the open problems and literature of this section, the reader may refer to \cite{LauraOb} and references therein. When studying the nilpotent cone, a few questions arise from the work of Gukov and his colleagues in relation to quantization:
\begin{open}
What is the brane quantization \`a la {\rm \cite{sergei1}} of the branes in the nilpotent cone mentioned in this section, and how does this relate to the curve quantization \`a la {\rm \cite{motohico}} of the spectral curves defined by Higgs bundles in those branes?
\end{open}

 \begin{open}
 Use the above methods to construct branes for wild Hitchin systems, and approach Langlands duality as appearing in~{\rm \cite{sergei2007}}.
 \end{open}

\subsection[Surface group representations and $GW$-components]{Surface group representations and $\boldsymbol{GW}$-components} Finally, it should be noted that branes which lie completely over the singular locus of the Hitchin fibration also play an important role in representation theory. In particular, the following has been predicted by Guichard and Wienhard:

\begin{Conjecture}[{\cite[Conjecture~5.6]{anna2}}]\label{A1} Additional connected components coming from positive representations $($through the notion of $\Theta$-positivity$)$, giving further families of higher Teichm\"uller spaces, appear in the moduli space of surface group representations into ${\rm SO}(p+q,p)$ for $q\geq1$.
\end{Conjecture}

In the case of the moduli space $\mathcal{M}_{{\rm SO}(p+1,p)}$, the existence of the extra Guichard--Wienhard components (or simply {\it GW-components}) as predicted in \cite[Conjecture~5.6]{anna2} is known to be true \cite{aparicio,brian2,brian3}, and moreover it was shown in~\cite{brian3} that those components indeed contain $\Theta$-positive representations.
 From the perspective of the spectral data description of the $(B,A,A)$-branes of ${\rm SO}(p+q,p)$-Higgs bundles of \cite{orthogonal1, orthogonal2}, natural candidates for $GW$-components for arbitrary $q\in \mathbb{N}$ are the following:

 \begin{Conjecture}[{\cite[Section 7]{orthogonal1}}]
The natural candidates for the $GW$-components of Conjecture~{\rm \ref{A1}} conjectured to exist by Guichard and Wienhard {\rm \cite[Conjecture~5.6]{anna2}} are those containing Higgs bundles whose spectral data\footnote{The spectral data is a triple $(L,M,\tau)$ consisting of a line bundle $L$, an orthogonal bundle $M$ (on an auxiliary curve) and an extension class~$\tau$.} $(L , M , \tau )$ in~{\rm \cite{orthogonal1}} has the form $( \mathcal{O} , \mathcal{O}^q , \tau )$. Alternatively, this is equivalent to taking ${\rm SO}(p+q,p)$-Higgs bundles whose vector bundle is of form $(W , V \oplus \mathcal{O}^{q-1})$, where the pair $(W,V)$ gives the vector bundles of one of the ${\rm SO}(p+1,p)$-Higgs bundles in the $GW$-components known to exist.
 \end{Conjecture}

To prove that this actually gives the $GW$-components, the monodromy action \`a la \cite{monodromy1, monodromy2} should be taken into consideration as well as the behavior over singular fibers. On the symplectic side, from the study of spectral data in \cite{orthogonal1}, one can see the geometric reason for the absence of any {\it extra} $GW$-components in the $(B,A,A)$-brane of ${\rm Sp}(2p+2q,2p)$-Higgs bundles.

\section{Higgs bundles and Calabi--Yau geometry} \label{sec-string}

Higgs bundles have played an important role in string theory in a wide range of contexts. But one recent application has provided some of the perhaps most surprising connections between previously unrelated aspects of geometry -- that is, the links between Hitchin systems and the geometry of singular $3$- and $4$-(complex)dimensional Calabi--Yau (CY) varieties.

\subsection{Calabi--Yau integrable systems and Hitchin Systems}The first hint of such a connection appeared in \cite{Diaconescu:2006ry} in which links were developed between Calabi--Yau integrable systems and Hitchin Systems. Briefly, as described in Section \ref{Higgs}, the Hitchin system forms an integrable system through the definition of the Hitchin fibration \eqref{hitchin_map}, whose generic fibers are even abelian varieties obtained through branched coverings of an underlying Riemann surface $\Sigma$ with genus $g\geq2$. On the other hand, \emph{Calabi--Yau} integrable systems were first explored for families of Calabi--Yau $3$-folds in \cite{donagimarkman, Donagi:1995am}, where the base of the system was formed by the moduli space of Calabi--Yau varieties in the family, and the fibers were formed by the Deligne cohomology groups of the intermediate Jacobians:
\begin{gather*}
J^2(X)=H^3(X, \mathbb{C})/\big(F^2H^3(X,\mathbb{C}) + H^3(X, \mathbb{Z})\big)
\end{gather*}
of the Calabi--Yau 3-folds $X$. Fiber and base fit together into a total space carrying a holomorphic symplectic form and the fibers are Lagrangian \cite{donagimarkman}. Furthermore, in remarkable work~\cite{Diaconescu:2006ry} Diaconescu, Donagi, and Pantev developed an isomorphism between Calabi--Yau integrable systems and those of Hitchin, the DDP correspondence. More precisely, by considering a~smooth projective complex curve $\Sigma$ and an ADE group~$G$, for a~fixed pair $(\Sigma, G)$ they constructed a~family of quasi-projective (i.e., \emph{non-compact}) CY $3$-folds (defined as $\operatorname{Tot}(V)$ for a rank 2 vector bundle, $V$, satisfying $\det (V)=K_{\Sigma}$).

Treating the moduli space of the non-compact CY manifold as the base of a Hitchin integrable system for the group $G$, a correspondence between the CY integrable system (whose fibers are the intermediate Jacobians of a family of non-compact CY 3-folds) and that of the Hitchin system (whose fibers are Prym varieties of the corresponding spectral covers) was explicitly laid out for the Lie groups $A_k$, but the description is only valid away from the discriminant. This mapping between Hitchin and CY integrable systems was nicely generalized via a sheaf-theoretic approach to the remaining simple Lie groups $B_k$, $C_k$, $F_4$ and $G_2$ in~\cite{beck}.

This important correspondence found a ready audience within string theory in the context of F-theory \cite{Vafa:1996xn} -- a geometric approach to compactifications of the type IIB string with non-trivial axio-dilaton backgrounds~-- in which the effective physics of the type IIB compactification to $(12-2n)$ spacetime dimensions is encoded in the geometry of an elliptically fibered (or more generally genus one fibered), complex Calabi--Yau $n$-fold, $\pi\colon X_n \to B_{n-1}$. The degeneration of the elliptic fibers encode information about a Lie group, $G$, corresponding to D7-branes wrapping the discriminant locus of the elliptic fibration (see Section~\ref{Mboyo_sec} for further details). From a physical perspective, the Calabi--Yau geometry is a tool to investigate intersecting brane theories, which are innately linked to Higgs bundles. Within a ``local'' description of F-theory, intersecting branes (wrapping a sub-variety $\Sigma \subset X_n$) come equipped with an adjoint field $\Phi$ (the Higgs field) which parametrizes normal motions of a stack of branes. Matter fields are fluctuations around the background $\langle \Phi \rangle$ and Yukawa couplings measure obstructions to extend these solutions beyond the linear order. Usually the Higgs field $\Phi$ is taken to live in a Cartan subspace of the Lie algebra so that only the eigenvalues of~$\Phi$ are relevant. But this seems to be an incomplete description in many situations relevant to interesting physical models.

\subsection{T-branes and Hitchin systems} Usually the connection with F-theory is made using the spectral cover, defined through $\{\det ( \Phi-\eta \operatorname{Id})=0\}$, and which in the case that $\Sigma \subset X$ (with $X$ a CY $n$-fold) locally defines the CY as the normal cone of $\Sigma$. However, the spectral cover will not accurately parameterize the local geometry of $X$ when the Higgs field is non-diagonalizable (for example when $\Phi$ is upper triangular)~\cite{Donagi:2003hh}. So-called \emph{T-branes} are non-Abelian bound states that generalize intersecting branes and admit a matrix of normal deformations (or Higgs field) that is nilpotent over some loci~\cite{Cecotti:2010bp,Donagi:2011jy,Donagi:2011dv}. Mathematically, T-branes correspond to singular fibers in the Hitchin fibration.

As they first originated in the physics literature, the Hitchin systems corresponding to T-branes were not explicitly linked geometrically to the background Calabi--Yau elliptic fibrations of F-theory. A first step in this direction was taken in \cite{Anderson:2013rka} which attempted to extend some of the links developed in \cite{Diaconescu:2006ry} to \emph{compact Calabi--Yau varieties}. A limiting mixed Hodge structure analysis was employed to study the form of the intermediate Jacobian of CY 3-folds in the limit that the geometry became singular. In certain singular limits of the elliptic fibration, the degeneration of the intermediate Jacobian, $J^2(X)$, leads to \emph{an emergent Hitchin system}~-- i.e., generates the moduli space of Higgs bundles defined over the discriminant locus of the fibration~\cite{Anderson:2013rka, Donagi:1995am}. In~\cite{Anderson:2013rka} T-branes were explored in the context of six-dimensional F-theory vacua, that is using compactifications of F-theory
on \emph{singular} elliptically fibered Calabi--Yau $3$-folds, $\pi\colon X_3 \to B_2$.

The intrinsic intersecting brane Hitchin system was defined over a curve $\Sigma \subset B$ in the base of the elliptic fibration, obtained through a component of the discriminant locus ($\Delta=0$) describing degenerating fibers. Upon a crepant resolution of the singular variety, it was argued that the geometric remnants of T-branes correspond to periods of the three-form potential of F-theory valued in the intermediate Jacobian of a now smooth Calabi--Yau 3-fold. Moreover, in~\cite{Anderson:2013rka} a~partial compactification of the DDP correspondence was established and it was demonstrated that the Hitchin system defined on the discriminant locus is contained in the local part of the (compact) Calabi--Yau integrable system:
\begin{gather}\label{cy_maps}\begin{split}
& \xymatrix@R=15pt@C=15pt{
&\pi^* \mathcal{M} \ar@{_{(}->}[dl] \ar[d] \ar[r] &
\mathcal{M} \ar[d]^{\rm Hit} &\\
M \ar[r] & \widetilde M_{\rm cx} \ar[r]^\pi & M_{\rm loc},}
\end{split}
\end{gather}
where $\mathcal{M}$ and $M$ are the full Hitchin and Calabi--Yau moduli spaces, respectively, and $\widetilde{M}_{\rm cx}$ and $M_{\rm loc}$ the complex
structure moduli spaces of the resolved Calabi--Yau geometry and local deformations of the singular Calabi--Yau variety (preserving the form of the singular elliptic fibers). The maps are defined such that the right map is the Hitchin fibration and the upper left (diagonal) map is an inclusion. This correspondence was established for singular CY 3-folds with $A_n$-type singularities. See \cite{Anderson:2013rka}, Section A.4 for details.

\begin{Conjecture}[\cite{Anderson:2013rka}] The correspondence described in~\eqref{cy_maps} can be extended to any compact, singular, elliptically fibered CY $3$-fold with singular fibers associated to ${\mathcal G}$-symmetry $($co-dimension~$1$ over the base and admitting a crepant resolution, see Section~{\rm \ref{Mboyo_sec}} for details$)$, and ${\mathcal H}$-type Hitchin system defined over the discriminant locus of the elliptic fibration such that ${\mathcal H} \subset {\mathcal G}$.\end{Conjecture}

Although we will not explore it in detail here, it is also expected that correspondences between Hitchin and CY moduli spaces (in either the compact or non-compact CY setting) should extend to the Deligne cohomology of Calabi--Yau $4$-fold geometries \cite{Bies:2014sra} (see also \cite{Bena:2017jhm,Collinucci:2016hpz,Collinucci:2017bwv,Collinucci:2014taa, Marchesano:2017kke} for recent progress on T-branes) and Higgs bundles defined over complex surfaces \cite{simpsonhodge}.

\subsection{Wild Hitchin systems and F-theory}An intrinsic difficulty with the Hitchin systems arising within F-theory comes from the fact that the Higgs bundles are defined on the discriminant locus of the elliptically fibered CY manifold~-- and hence not on smooth Riemann surfaces, but rather on complex curves that are in general singular (including sometimes non-reduced and reducible). At the singular/intersection points of such a curve, the physical theory suggests that the associated Higgs bundles should also exhibit singularities. That is, in this context, it is natural to also consider stable pairs with \emph{singular connections} (see also the discussion in Section~\ref{sec-curves}).

Thus far, work has focused primarily on so-called parabolic Hitchin systems \cite{konno93, simpsonhodge, simpsonhiggs} which accommodate the possibility of simple poles in the gauge and Higgs fields at marked points on a Riemann surface. However, many questions~-- of both mathematical as well as physical interest~-- require the consideration of higher order singularities in the gauge fields. Wild/irregular Higgs bundles \cite{AkerSzabo,2012arXiv1203.6607B,1996alg.geom.10014B, FredricksonNeitzke} extend this formalism to include stable, integrable connections with irregular singularities of the form
\begin{gather*}
d+A_n \frac{dz}{z^n}+ \cdots +A_1 \frac{dz}{z},
\end{gather*}
with $n>1$ and stable parabolic Higgs pairs $(E,\Phi)$ where the Higgs field has polar parts, e.g.,
\begin{gather*}
T_n \frac{dz}{z^n}+ \cdots + T_1 \frac{dz}{z}.
\end{gather*}
As in Simpson's construction \cite{simpsonhiggs} for parabolic Higgs bundles, a natural assumption for this study is that the connections and Higgs fields are holomorphically gauge equivalent to ones with diagonal polar parts (this is weakened slightly in \cite{FredricksonNeitzke}). This leads to a correspondence between singularities (after diagonalizing) of the form $T_i=\frac{1}{2}A_i$ for $i\geq 2$~\cite{biquardboalch}. The moduli space of such wild Higgs bundles was described in \cite{2012arXiv1203.6607B,boalch2} as a hyperk\"ahler quotient. Already these irregular Hitchin systems have played a~significant role in the geometric Langlands program \cite{2006math.....11294F,sergei2007}, and string applications including topologically twisted $\mathcal{N}=4$ super Yang--Mills theories \cite{sergei2007,kapustinwitten}, particularly in so-called 'Stokes phenomena' \cite{Witten:2007td,Xie:2012hs} (which describe how the asymptotical behavior of the solutions changes in different angular regions around the singularity. Stokes matrices link the solutions in different regions and define a generalized monodromy which plays a central role in describing a wild Hitchin moduli space).

Recent progress \cite{Anderson:2017rpr} has demonstrated that the study of ordinary smooth Hitchin systems is insufficient in the context of $6$-dimensional F-theory compactifications. Not only should generic CY $3$-folds have a correspondence to parabolic or wild Hitchin systems, but in general deformations of the singular variety ($M_{\rm loc}$ in \eqref{cy_maps} above) can {\it dynamically change the pole order} of the relevant singular complexified connections appearing in the Higgs bundles (see \cite{Anderson:2015cqy} and \cite[Section~5.1]{Anderson:2017rpr} for an example). That is, by varying the complex structure moduli of a singular CY variety, the location of singularities in the discriminant locus $\Delta$ can be tuned to coincide. This tuning, when viewed from the intersecting brane models, should correspond to a parametric deformation of a Hitchin System in which the location of simple poles are tuned and forced to coincide into higher order poles. In many instances it seems this tuning can be done without changing the dimension of the underlying CY/Hitchin moduli space.

\begin{Conjecture}[\cite{Anderson:2017rpr}] There exists a flat morphism between the moduli spaces ${\mathcal M}_{\rm par}$ and ${\mathcal M}_{\rm wild}$ in the case of a singular parabolic Higgs bundle with $n$-simple poles in its connection and that of a wild Higgs bundle with a single, higher order pole of order $n$.\end{Conjecture}

This limiting process in the context of CY varieties also leads to the open question:

\begin{open}[\cite{Anderson:2017rpr}] Does the Stokes phenomenon exhibited by wild Higgs bundles have an analog in the moduli space of singular, elliptically fibered CY geometries or CY integrable systems?\end{open}

\subsection{Singular CY varieties}Finally, it should be noted that there are likely many unexplored links between classification problems in parabolic/wild Hitchin Systems and singular CY varieties. In general, the criteria for CY 3-folds to exhibit a generic singularity everywhere in its complex structure moduli space has attracted interest from the physics community in the context of ``non-Higgsable clusters'' \cite{Anderson:2014gla,Halverson:2015jua,Morrison:2012np,Morrison:2014lca}. The maps in~\eqref{cy_maps} embedding the Hitchin moduli space into that of the singular CY variety indicate that in such cases the highly constrained form of the singular CY geometry must correspond to an equally constrained Hitchin system. From the underlying effective physics of F-theory, this correspondence is linked to an ${\rm SU}(2)$ R-symmetry which can rotate components of hypermultiplets of the $6$-dimensional effective theory \cite{Aspinwall:1998he}. Here these halves of hypermultiplets correspond to complex structure moduli and degrees of freedom in the intermediate Jacobian of the CY variety (the so-called ``RR-moduli''), respectively. Thus, any T-brane solution (or more generally Higgs bundle on the brane with associated spectral cover) must correspond under hypermultiplet rotation to a deformation of complex structure of the singular 3-fold \cite{Anderson:2013rka}. In particular, in the case of non-Higgsable CY $3$-fold geometries an open question is to understand the following:

\begin{open}In the case that a {\it non-Higgsable} CY manifold exhibits ${\mathcal G}$ singular elliptic fibers over a $($possibly singular$)$ curve $\Sigma$, does this correspond to a trivial/empty Hitchin moduli space $($including parabolic or wild Hitchin moduli space depending on the singularities of $\Sigma)$ of ${\mathcal H}$-Higgs bundles over $\Sigma$ where ${\mathcal H} \subset {\mathcal G}$?\end{open}

One example where this open question can be confirmed in the affirmative is in the case of CY $3$-folds $\pi\colon X \to \mathbb{F}_n$, defined as a Weierstrass model (see Section~\ref{Mboyo_sec}) over a Hirzebruch surface, $\mathbb{F}_n$. For each $n>2$, such elliptically fibered geometries are generically singular \cite{Bershadsky:1996nh, Morrison:1996pp}. For $n=3$ for example there is a generic ${\rm SU}(3)$ singularity (more specifically Kodaira type IV fibers) over a discriminant locus $\Delta \subset \mathbb{F}_3$ which takes the form of a smooth curve of genus zero (one of the sections of the rational fibration of~$\mathbb{F}_3$). The fact that this symmetry is ``un-Higgsable'' in the physical theory corresponds to the triviality of ${\rm SL}(N, \mathbb{C})$-Higgs bundles over $\mathbb{P}^1$ with $N=2,3$. See~\cite{Anderson:2013rka} for further examples involving parabolic Hitchin systems over $\mathbb{P}^1$ with marked points.

In general, there are a large number of possible connections between Higgs bundles and the effective theories and geometry arising within string/M-/F-theory. One new correspondence has recently arisen within the context of $4$-dimensional compactifications of F-theory in ${\mathcal N}=4$ supersymmetric Yang--Mills theories (with unity gauge groups) which are quotiented by particular combinations of $R$-symmetry and ${\rm SL}(2, \mathbb{Z})$ automorphisms (such theories can arise as D3-branes probing terminal singularities in F-theory) \cite{Aharony:2015oyb,Aharony:2016kai,Arras:2016evy, Garcia-Etxebarria:2015wns,Garcia-Etxebarria:2016erx}.

 \begin{open}What generalizations of Higgs bundles correspond to the new ${\mathcal N}=3$ supersymmetric theories recently discovered in~{\rm \cite{Garcia-Etxebarria:2015wns}}?
\end{open}

The self duality equations of ${\mathcal N}=4$ supersymmetric Yang--Mills theories have led to a rich interplay between theories of branes arising in string theory and Higgs bundles. It would be intriguing to understand whether such links could arise between ${\mathcal N}=3$ theories and ``cousins'' of the Hitchin system over Riemann surfaces. Finally, it should be noted that within F-theory and the subject of T-branes there remain many open questions linking so-called "matrix factorization" techniques, K-theory, Hitchin systems and Calabi--Yau geometry (see, e.g.,~\cite{Bena:2017jhm}).

\section[Elliptic fibrations, Weierstrass models, and Calabi--Yau resolutions]{Elliptic fibrations, Weierstrass models,\\ and Calabi--Yau resolutions}\label{Mboyo_sec}

 Since elliptic fibrations play an important role when studying the relations between Higgs bundles and F-theory, we shall conclude these notes with a review of some of the basic ideas and recent advances. The reader should not take this a thorough review, but rather a brief, curated overview of some essential aspects of the underlying geometry.

\begin{Definition} A surjective proper morphism $\varphi\colon Y\to B$ between two algebraic varieties $Y$ and $B$ is called an {\it elliptic fibration} if the generic fiber of $\varphi$ is a smooth projective curve of genus one and $\varphi$ has a rational section. When $B$ is a curve, $Y$ is called an {\it elliptic surface}. Moreover, when $B$ is a surface, $Y$ is said to be an {\it elliptic $3$-fold}. In general, if $B$ has dimension $n-1$, $Y$~is called an elliptic $n$-fold.
\end{Definition}

\subsection{Classification of singular fibers}The locus of singular fibers of an elliptic fibration, $\varphi\colon Y\to B$, is called the {\it discriminant locus}, and is denoted by $\Delta(\varphi)$, or simply $\Delta$ when the context is clear. If the base~$B$ is smooth, the discriminant locus is a divisor \cite{Dolgavcev.Purity}. In the early 1960s, Kodaira classified singular fibers of minimal elliptic surfaces in terms of numerical invariants showing that there are 8 possibilities including two infinite series and 6 exceptional cases~\cite{Kodaira,Kodaira+}. Soon after, N\'eron obtained an equivalent classification in an arithmetic setting using explicit regularizations of singular Weierstrass models~\cite{Neron}. Based on N\'eron's analysis, Tate proposed an algorithm that allows (among other things) the determination of the type of singular fibers of a Weierstrass model by analyzing the valuation of its coefficients~\cite{Tate}. Under appropriate conditions, Kodaira's classification of singular fibers of an elliptic surface and Tate's algorithm can be used to describe the possible singular fibers and monodromies of an elliptically fibered $n$-fold over points in codimension-1 in the base.

However, over points in codimension-2, new fibers not in Kodaira's list are known to occur \cite{Miranda,Szydlo}. These are frequently referred to as {\em collisions of singular fibers} as they usually appear at the intersections of two divisors of the discriminant locus of the elliptic fibration. In general however, no classification exists for the singular fibers of elliptic $3$-folds and $4$-folds, leading to the broad question:

\begin{open}How can one geometrically classify non-Kodaira fibers?
\end{open}

Under some assumptions, the answer is known for the elliptic $n$-folds called Miranda models \cite{Miranda, Szydlo}. More generally, in the case of flat elliptic fibrations obtained by crepant resolutions of Weierstrass models, non-Kodaira fibers are expected to be contractions of usual Kodaira fibers. This is proven for elliptic 3-folds~\cite{Cattaneo:2013vda} and confirmed in all known examples of non-Kodaira fibers appearing in the F-theory literature \cite{EFY,EKY,ESY1,ESY2,EY,Lawrie:2012gg, Morrison:2011mb,Tatar:2012tm}.

 In view of the links described in Section~\ref{sec-string} between elliptic CY geometry and Hitchin systems, this leads naturally to the conjecture that two classification problems might be linked:

\begin{open}How is the classification of non-Kodaira fibers of an elliptically fibered Calabi--Yau $3$-fold related to the classification of parabolic or wild Hitchin systems defined over the discriminant locus?
\end{open}

\subsection{Weierstrass models} Since an elliptic fibration over a smooth base is birational to a (possibly singular) {\it Weierstrass model}~\cite{Deligne.Formulaire}, the starting point of such an analysis will usually be a Weierstrass model. We shall review here the main features of these models, following the notation of Deligne~\cite{Deligne.Formulaire}. Let $\mathcal{L}$ be a line bundle over a quasi-projective variety~$B$. We define the following projective bundle (of lines):
\begin{gather*}
\pi\colon \ X_0=\mathbb{P}_B\big[\mathcal{O}_B\oplus \mathcal{L}^{\otimes 2}\oplus \mathcal{L}^{\otimes 3}\big]\longrightarrow B.
\end{gather*}
We denote by $\mathcal{O}_{X_0}(1)$ the dual of the tautological line bundle of the projective bundle $X_0$. The relative projective coordinates of $X_0$ over $B$ are denoted $[z:x:y]$, where $z$, $x$, and $y$ are defined respectively by the natural injection of $\mathcal{O}_B$, $\mathcal{L}^{\otimes 2}$, and $\mathcal{L}^{\otimes 3}$ into $\mathcal{O}_B\oplus \mathcal{L}^{\otimes 2}\oplus \mathcal{L}^{\otimes 3}$. Hence, $z$ is a~section of $\mathcal{O}_{X_0}(1)$, $x$ is a section of $\mathcal{O}_{X_0}(1)\otimes \pi^\ast \mathcal{L}^{\otimes 2}$, and~$y$ is a section of $\mathcal{O}_{X_0}(1)\otimes \pi^\ast \mathcal{L}^{\otimes 3}$. The most general Weierstrass equation is then the zero locus of the following section of $\mathcal{O}(3)\otimes \pi^\ast \mathcal{L}^{\otimes 6}$ in $X_0$
\begin{gather*}
F=y^2z+ a_1 xy z + a_3 yz^2 -\big(x^3+ a_2 x^2 z + a_4 x z^2 + a_6 z^3\big),
\end{gather*}
where $a_i$ is a section of $\pi^\ast \mathcal{L}^{\otimes i}$. The line bundle $\mathcal{L}$ is called the {\em fundamental line bundle} of the Weierstrass model $\varphi\colon Y\to B$ and can be defined directly from the elliptic fibration $Y$ as $\mathcal{L}=R^1 \varphi_\ast \mathcal{O}_Y$. The Weierstrass model has a trivial canonical class when the fundamental line bundle $\mathcal{L}$ is the anti-canonical line bundle of $B$.

Each crepant resolution of a singular Weierstrass model is a relative minimal model (in the sense of the Minimal Model Program) over the Weierstrass model~\cite{Matsuki.Weyl}. When the base of the fibration is a curve, the Weierstrass model has a unique crepant resolution. On the other hand, when the base is of dimension two or higher, a crepant resolution does not always exist; furthermore, when it does, it is not necessarily unique. Different crepant resolutions of the same Weierstrass model are connected by a finite sequence of flops (see for example \cite{EJaK,EJeK,ESY1,ESY2,EY,Hayashi:2014kca, Krause:2011xj, Matsuki.Weyl}). Crepant resolutions of Weierstrass models have the same Euler characteristic, and these have recently been computed in \cite{Euler}.

Following F-theory, we can attach to a given elliptic fibration a Lie algebra~$\mathfrak{g}$, a representation~$\mathbf{R}$ of~$\mathfrak{g}$, and a hyperplane arrangement ${\rm I}(\mathfrak{g},\mathbf{R})$. The Lie algebra~$\mathfrak{g}$ and the representation $\mathbf{R}$ are determined by the fibers over codimension-1 and codimension-2 points, respectively, of the base in the discriminant locus. The hyperplane arrangement ${\rm I}(\mathfrak{g},\mathbf{R})$ is defined inside the dual fundamental Weyl chamber of $\mathfrak{g}$ (i.e., the dual cone of the fundamental Weyl chamber of~$\mathfrak{g}$), and its hyperplanes are the set of kernels of the weights of $\mathbf{R}$. Moreover, one may study the network of flops using the hyperplane arrangement ${\rm I}(\mathfrak{g}, \mathbf{R})$ inspired from the theory of Coulomb branches of five-dimensional supersymmetric gauge theories with eight supercharges~\cite{IMS}.

The network of crepant resolutions is isomorphic to the network of chambers of the hyperplane arrangement I$(\mathfrak{g}, \mathbf{R})$ defined by splitting the dual fundamental Weyl chamber of the Lie algebra $\mathfrak{g}$ by the hyperplanes dual to the weights of $\mathbf{R}$. The hyperplane arrangement I$(\mathfrak{g}, \mathbf{R})$, its relation to the Coulomb branches of supersymmetric gauge theories and the network of crepant resolutions are studied, among others, in
 \cite{Braun:2014kla,Diaconescu:1998cn,EJJN1,EJJN2,EJaK,EJeK,ESY1,ESY2,Grimm:2011fx,Hayashi:2014kca,Hayashi:2013lra, IMS}.

The representation $\mathbf{R}$ attached to an elliptic fibration can be derived systematically using intersection theory \cite{Aspinwall:1996nk}.
Indeed, let $C$ be a vertical curve, i.e., a curve contained in a fiber of the elliptic fibration. Let $S$ be an irreducible component of the reduced discriminant of the elliptic fibration $\varphi\colon Y\to B$. The pullback of $\varphi^* S$ has irreducible components $D_0, D_1, \ldots, D_n$, where $D_0$ is the component touching the section of the elliptic fibration. The {\em weight vector} of~$C$ over~$S$ is by definition the vector ${\varpi}_S(C)=(-D_1\cdot C, \ldots, -D_n\cdot C)$ of intersection numbers $D_i\cdot C$ for $i=1,\ldots, n$. To an elliptic fibration, we associate a representation $\mathbf{R}$ of the Lie algebra $\mathfrak{g}$ as follows. The weight vectors of the irreducible vertical rational curves of the fibers over codimension-2 points form a set $\Pi$ derived by intersection theory. The saturation of $\Pi$ (by adding and subtracting roots) defines uniquely a representation $\mathbf{R}$. This method due to Aspinwall and Gross \cite[Section~4]{Aspinwall:1996nk} explains how the representation $\mathbf{R}$ can be deduced even in presence of non-Kodaira fibers~\cite{Marsano:2011hv}. The method can be formalized using the notion of {\em saturation set of weights} borrowed from Bourbaki~\cite{EJaK,EJeK}.

One interesting property of the derivation of the representation $\mathbf{R}$ from intersection theory is that it does not assume the Calabi--Yau condition nor relies on anomaly cancellations. Hence, from that point of view, the representation attached to an elliptic fibration is
purely a geometric data of the elliptic fibration that also controls aspects of its birational geometry via the hyperplane arrangement ${\rm I}(\mathfrak{g},\mathbf{R})$. There are subtleties in presence of exotic matter \cite{Anderson:2015cqy}, when the component of the discriminant supporting the gauge group is singular \cite{Klevers:2016jsz}, in presence of a~non-trivial Mordell--Weil group \cite{Mayrhofer:2014opa}, when the codimension-two fibers are non-split \cite{EJaK,EJeK}, or when the fibration is non-flat \cite{Lawrie:2012gg}. Although we understand the structure of the hyperplane arrangement ${\rm I}(\mathfrak{g},\mathbf{R})$ for most of the F-theory models with simple groups (see for example \cite{Diaconescu:1998cn,EJJN1,EJJN2,EJaK,EJeK,Hayashi:2014kca, IMS}), the structure in presence of semi-simple groups is still not well explored. This lead to the following question.

\begin{open}What are the intersection properties of $($exotic$)$ representations appearing in F-theory and the structure of their associated hyperplane arrangements?
\end{open}

\subsection{Superconformal field theories in the context of F-theory}Finally, it should be noted that there is a rich array of open questions that have arisen form the recent investigations into superconformal field theories in the context of F-theory. The superconformal algebra is a graded Lie algebra that combines the conformal Poincar\'e algebra and supersymmetry. Some of the most basic data characterizing a superconformal field theory (SCFT) is the number of spacetime dimensions in which it is defined, and the amount of supersymmetry generators and their chirality.

Recently, substantial interest has centered on six-dimensional SCFTs with $(2,0)$ and $(1,0)$ supersymmetry. According to the seminal work of Werner Nahm, SCFTs are only possible for spacetime dimensions~2, 3, 4, 5 and~6~\cite{Nahm:1977tg}. In particular, the (2,0) theories in $d = 6$ are the SCFTs with the maximal amount of supersymmetry in the highest dimension~\cite{Tachikawa:2013kta}. The six-dimensional superconformal field theories with (1,0) supersymmetry are among the least understood quantum field theories, for example, they do not always have a Lagrangian formulation \cite{Seiberg:1996vs,Witten:2002eh}. They are connected to questions in broad areas such as Donaldson-Thomas theory of Calabi--Yau manifolds, modular and automorphic forms \cite{Gu:2017ccq, Haghighat:2014vxa}, singularities \cite{DelZotto:2014hpa,Heckman:2015bfa, Heckman:2013pva}, quivers \cite{Aganagic:2016jmx,Bhardwaj:2015xxa,Kimura:2016dys}, and representation theory \cite{Aganagic:2017smx}. As they arise in F-theory and CY elliptic fibrations, it is then natural to ask:

\begin{open}What is the geometry of elliptic fibrations used to model $(1,0)$ theories? How are the conformal matter connected to the structure of Higgs bundles appearing in F-theory?
\end{open}

The crepant resolutions of the singularities of CY elliptic fibrations exhibiting SCFT loci provide a beautiful connection between the mechanism of anomaly cancellation as seen in physics and topological quantities that have been recently discovered. An understanding of the SCFT geometry must be linked to the simplest building blocks of (1,0) theories, the so called non-Higgsable clusters \cite{Heckman:2015bfa,Heckman:2013pva, Morrison:2012np,Morrison:2012js}. First steps towards the analysis of the crepant resolutions of such SCFT loci (including so-called "matter transitions" in F-theory \cite{Anderson:2015cqy}) is already underway \cite{DelZotto:2017pti,EJaK,EJeK, Haghighat:2014vxa,Lawrie:2012gg}.

\subsection*{Acknowledgements}
The authors would like to thank the American Institute of Mathematics for the support and hospitality which made the 2017 workshop on {\it Singular Geometry and Higgs Bundles in String Theory} possible (and on which this survey of ideas/open questions is based). In addition we would like to thank Steven Rayan for helpful comments on the manuscript. The work of L.B.~Anderson is supported in part by NSF grant PHY-1720321 and is part of the working group activities of the 4-VA initiative ``A Synthesis of Two Approaches to String Phenomenology''. M.~Esole is supported in part by the National Science Foundation (NSF) grant DMS-1701635 ``Elliptic Fibrations and String Theory''. The work of L.P.~Schaposnik is partially supported by the NSF grant DMS-1509693, and by the Alexander von Humboldt Foundation.

\pdfbookmark[1]{References}{ref}
\LastPageEnding

\end{document}